\documentclass[11pt]{amsart}
\usepackage{amsmath,amsfonts,amsthm,amscd,amssymb,mathrsfs,amssymb,bm}
\usepackage{mathrsfs}
\usepackage[dvips]{graphicx}
\usepackage{wrapfig}
\usepackage[all]{xy}
\newtheorem{theorem}{Theorem}

\newtheorem{proposition}[theorem]{Proposition}
\newtheorem{lemma}[theorem]{Lemma}

\newtheorem{remark}[theorem]{Remark}
\newcommand{\CP}{\mathbb{CP}}
\newcommand{\CC}{\mathbb{C}}

\newcommand{\ol}{\overline}
\newcommand{\lra}{\longrightarrow}

\newcommand{\proofend}{\hfill$\square$}

\newcommand{\Bs}{{\rm{Bs}}}

\newcommand{\vsp}{\vspace{3mm}}

\numberwithin{equation}{section}
\numberwithin{theorem}{section}
\begin{document}
\bibliographystyle{alpha} 
\title{Classification of Moishezon twistor spaces on $4\mathbb{CP}^2$}
\author{Nobuhiro Honda}
\address{Mathematical Institute, Tohoku University,
Sendai, Miyagi, Japan}
\email{honda@math.tohoku.ac.jp}
\thanks{The author was partially supported by the Grant-in-Aid for Young Scientists  (B), The Ministry of Education, Culture, Sports, Science and Technology, Japan. }
\begin{abstract}
In this paper we provide a classification of all Moishezon twistor spaces
on the connected sum of four complex projective planes.
This is given by means of the anticanonical system
of the twistor spaces.
In particular, we show that the anticanonical map is birational, two to one over the image,
or otherwise the image of the anticanonical map is a rational surface.
We also obtain the structure of the images of the anticanonical map 
in each of the three cases in  quite concrete forms.
\end{abstract}
\maketitle
\setcounter{tocdepth}{1}
\vspace{-5mm}

\section{Introduction}
By Taubes' theorem \cite{Ta} there exist a huge number of compact oriented 4-manifolds
which admit a self-dual conformal structure.
Associated to any self-dual structure is a 3-dimensional complex manifold of a special kind,
which is so called the {\em twistor space},
but
by a theorem of Hitchin \cite{Hi81} a  compact twistor space does not admit  a K\"ahler metric except
 two 
well-known examples.
Also, by a theorem of Campana \cite{C91},  the twistor space can be Moishezon only when the 4-manifold
is $S^4$ or $n\CP^2$, the connected sum of $n$ copies of complex projective planes.
By Kuiper \cite{Kui}, the standard metric on $S^4$ is the unique self-dual  structure on $S^4$.
Similarly, by Poon \cite{P86}, the Fubini-Study metric on $\CP^2$ is the unique self-dual structure on $\CP^2$ whose twistor space is Moishezon (or, of positive scalar 
curvature).
On $2\CP^2$, Poon \cite{P86} constructed a family of Moishezon twistor spaces parametrized by an open interval in $\mathbb  R$, and showed that if a self-dual structure is of positive scalar curvature, then 
the twistor space must be one of the twistor spaces in this family.
Thus the classification of Moishezon twistor spaces is completely over up to $2\CP^2$.

For the case $3\CP^2$, by the works of Poon \cite{P92} and Kreussler-Kurke \cite{KK92} a classification is given by means of the complete linear system $|K^{-1/2}|$, where
$K^{-1/2}$ is the natural square root of
the anticanonical line bundle, which is available on any twistor space.
Namely Moishezon twistor spaces on $3\CP^2$ can be classified into two types according to  whether 
$|K^{-1/2}|$ is base point free;
if free, the associated morphism is a degree 2 morphism onto $\CP^3$ and the branch divisor is a quartic surface,
whereas if not free, the image of the rational map is a non-singular quadratic surface in $\CP^3$,
and the twistor space has to be a LeBrun twistor space constructed in  \cite{LB91}
whose structure is also well understood.

For general $n\CP^2$, although a lot of Moishezon twistor spaces are known,
their classification seems still difficult.
In this paper we give a classification of Moishezon twistor spaces on $4\CP^2$, by means of the anticanonical system
of the twistor spaces.
We note that on $4\CP^2$, the system $|K^{-1/2}|$ is not enough for analyzing structure of twistor spaces, since in most cases the system  is  just a pencil.
A simplified form of our classification can be presented as follows:

\begin{theorem}\label{thm:main1}
Let $Z$ be a Moishezon twistor space on $4\CP^2$
and $\Phi$ the anticanonical map of $Z$.
Then exactly one of the following  three situations occurs:
\begin{enumerate}
\item[(I)] $\Phi$ is birational over the image,
\item[(II)] $\Phi$ is (rationally) two to one over a scroll of 2-planes over a conic
(so it is in $\mathbb{CP}^4$),
and the branch divisor restricts to a quartic curve on a general 2-plane of the scroll,
\item[
(III)] the image $\Phi(Z)$ is 2-dimensional, and general fibers of \,$\Phi$ \!are non-singular rational curves.
\end{enumerate}
\end{theorem}

Thus the anticanonical map nicely describes the structure of any Moishezon 
twistor spaces on $4\CP^2$.
Following these structure, let us call  a Moishezon twistor space $Z$ on $4\CP^2$ to be 
{\em double solid type} if $Z$ belongs to the case (II), and 
 {\em conic bundle type} if $Z$ belongs to the case (III).
(For the case (I) we simply say that the anticanonical map is birational.)
Then what we actually do is to investigate each of these three cases much more in detail and to clarify the structure of the anticanonical map $\Phi$ as follows.

For twistor spaces on $4\CP^2$ whose anticanonical map is birational, we show that the dimension of the anticanonical system is either 8 or 6.
We also show that if the dimension is 8 the anticanonical  image
is a (non-complete) intersection of 10 quadratic hypersurfaces in $\CP^8$, whose degree is 12,
and that if the dimension is 6, the anticanonical image is a complete intersection of three quadratic
hypersurfaces in $\CP^6$ (Theorem \ref{thm:birational}).

For twistor spaces of the double solid type, we classify them into 4 subtypes, 
according to the number of  irreducible components of the base curve 
of the system $|K^{-1/2}|$ (which will turn our to be a pencil). 
Distinction for these 4 kinds of spaces is significant since the defining equation 
of the branch divisor of the anticanonical map takes different forms for each cases.
But since determination of the equation requires much more detailed  analysis for the spaces,
we will discuss this topic in a different paper.

For twistor spaces of conic bundle types, 
we prove that the dimension of the anticanonical system is 8, 5, or 4.
We also show that if the system is 8-dimensional the anticanonical image is an embedded image of a non-singular quadric surface in $\mathbb{CP}^3$ by the Veronese embedding $
\CP^3\subset\CP^9$ induced from $|\mathscr O(2)|$,
and that if 5-dimensional the anticanonical image is the Veronese surface $\CP^2$ in $\CP^5$.
We also prove that if the anticanonical system is 4-dimensional the anticanonical image is an intersection of 2 hyperquadrics in $\CP^4$, and determine the quadratic defining polynomials in explicit forms
(Theorem \ref{thm:conicbundle}).
Thus the structure of the anticanonical image is completely determined.

We also discuss  where various known Moishezon twistor spaces on $4\CP^2$ are placed in this 
classification.
As a consequence we find that there are some new Moishezon twistor spaces,
but they do not occupy  a large part of the moduli space of all Moishezon twistor spaces.

\vsp
\noindent
{\bf Notations.}
The natural square root of the anticanonical bundle $K^{-1/2}$ on a twistor space is 
denoted by $F$.
This is called the fundamental line bundle.
For a line bundle $\mathscr L$ on a compact complex manifold, we write $h^i(\mathscr L)= \dim H^i(\mathscr L)$.
The dimension of $|\mathscr L|$ always refers $h^0(\mathscr L) -1 $.
${\rm{Bs}} |\mathscr L|$ denotes the base locus of the complete linear system $|\mathscr L|$.
For a non-zero element $s\in H^0(\mathscr L)$, $(s)$ means the zero divisor of $s$.
A curve on $\CP^1\times\CP^1$ of bidegree $(a,b)$ is called
 an $(a,b)$-curve.

\section{Classification of the twistor spaces on $4\CP^2$ by the anticanonical system}
\subsection{The base locus of the fundamental system with $h^0(F)=2$.}\label{ss:bc}
Let $Z$ be a (not necessarily Moishezon) twistor space on $4\mathbb{CP}^2$.
We suppose that the corresponding self-dual structure on $4\CP^2$ is of positive type in the sense that
 the scalar curvature is positive.
 (Note that by a theorem of Poon \cite{P88} this positivity always holds if $Z$ is Moishezon.)
Then thanks to a Hitchin's vanishing theorem \cite{Hi80}, the Riemann-Roch formula for the
line bundle $F$ gives
\begin{align}\label{RR1}
h^0(F) - h^1(F) = 2.
\end{align}
In particular, we always have $h^0(F)\ge 2$.
The structure of $Z$ satisfying $h^0(F)>2$ is well understood by the following result due to Kreussler \cite{Kr98}:
\begin{proposition}\label{prop:Kr}
 (i)
If $h^0(F)\ge 4$, then $h^0(F)= 4$ and $Z$ is a LeBrun twistor space (which is Moishezon).
(ii) If $h^0(F)=3$, $Z$ is Moishezon if and only if $\Bs\,|F|\neq\emptyset$.
Moreover, in this situation, $Z$ has to be a twistor space studied by Campana and Kreussler \cite{CK98} (which is also 
Moishezon).
\end{proposition}

Because generic twistor spaces on $4\mathbb{CP}^2$
satisfies $h^0(F)=2$, the twistor spaces in Proposition \ref{prop:Kr}
 are  quite special ones  among all Moishezon twistor spaces
on $4\CP^2$..
In order to classify Moishezon twistor spaces on $4\CP^2$, we have to investigate those which satisfy
$h^0(F)=2$.
Let $Z$ be such a twistor space.
Then because $Z$ is not a LeBrun twistor space, $Z$ does not have a pencil of degree one divisors by \cite{P92}.
So general members of the pencil $|F|$ is irreducible.
Take any real irreducible member $S$ of this pencil.
Then by a theorem of Pedersen-Poon \cite{PP94}, $S$ is a non-singular rational surface satisfying 
$K_S^2=0$.
We show the following.
(See also \cite{Kr98}.)

\begin{proposition}
\label{prop:cycle}
Let $Z$ be a Moishezon twistor space on $4\mathbb{CP}^2$ satisfying $h^0(F)=2$,
and  $S$  any  real irreducible member of the pencil $|F|$ as above.
Then $h^0(S,K_S^{-1}) = 1$.
Let $C$ be the unique anticanonical curve on $S$.
Then $\Bs\,|F| = C$ and $C$ is a cycle of non-singular rational curves. 
Moreover, the number $m$ of irreducible components of $C$ is even with $4\le m\le 12$.
\end{proposition}

Here, by a cycle of non-singular rational curves, we mean a connected reduced 
divisor on a non-singular surface, which is of the form
$$ \sum_{i=1}^m C_i \,\,{\text{ with }}\, m>2$$
where $C_1,\dots,C_m$ are mutually distinct smooth rational curves satisfying $C_i\cdot C_{i+1}=1$ for any $1\le i\le m$
while $C_{m+1} =C_1$.
(Actual situations we will encounter are the case where $m$ is even with $\ge 4$.)

\noindent {\em Proof of Proposition \ref{prop:cycle}.}
By \eqref{RR1} we have $H^1(F)=0$.
So the exact sequence 
$ 0 \to \mathscr O_Z\stackrel{\otimes s}{\to} F \to F|_S\simeq K_S^{-1} \to 0$,
where $s$ being a section of $F$ satisfying $(s)=S$, gives an exact sequence
\begin{align}
 0 \lra \mathbb C \stackrel{s\otimes}{\lra} H^0(F) \lra H^0(K_S^{-1}) \lra 0,
\end{align}
which means $\Bs\,|F| = C$.
Let $\epsilon:S\to \mathbb {CP}^1\times\mathbb{CP}^1$ be 
a birational morphism preserving the real structure which maps twistor lines in $S$ to   $(1,0)$-curves. 
Then since the image $\epsilon_*(C)$ (= the image as a divisor) is also an anticanonical curve, 
it is a $(2,2)$-curve, which is real.
Since twistor lines do not have a real point, this implies that if $\epsilon_*(C)$ has a non-reduced component,
then it must be a twice of a real $(1,0)$-curve.
But this contradicts $h^0(K_S^{-1})=1$, since $\epsilon$ cannot blowup points on the images of twistor lines
as $N_{L/S}\simeq\mathscr O_C$ for a twistor line $L$ on $S$.
Therefore $\epsilon_*(C) = \epsilon(C)$ and it must be of the form
\begin{align}\label{001}
C_1 + \ol C_1, {\text{ where $C_1$ is a non-real irreducible $(1,1)$-curve}},
\end{align}
\begin{align}\label{002}
C_1+C_2,  {\text{where $C_1$ is a real $(1,2)$-curve, and $C_2$ is a real $(1,0)$-curve}},
\end{align}
or
\begin{align}\label{003}
C_1 + C_2 + \ol C_1 + \ol C_2,  {\text{ where $C_1$ is a non-real $(0,1)$-curve, and $C_2$ is a non-real $(1,0)$-curve.}}
\end{align}
But again as $h^0(K_S^{-1}) = 1$, the  situation \eqref{002} cannot happen.
For the  situation \eqref{003}  it is obvious that $\epsilon(C)$ is a cycle of smooth rational curves.
The same conclusion holds for the case \eqref{001} since $C_1$ and $\ol C_1$ cannot touch at a point
since the induced real structure on $\CP^1\times\CP^1$ does not have a real point. 
Hence $\epsilon(C)$ is a cycle of smooth rational curves.
Moreover since $h^0(K_S^{-1})=1$, each of the single blowup
of the birational morphism $\epsilon:S\to\CP^1\times\CP^1$ always blows up a point on (the inverse image of) $\epsilon(C)$.
This means that $C$ is also a cycle of smooth rational curves and  that the number $m$ of irreducible components
satisfies $2\le m\le 4 + 8 = 12$.
Obviously $m$ is even by the real structure.
It remains to see $m\neq 2$.
If $m=2$, every step of $\epsilon$ blows up a smooth point of the anticanonical cycle \eqref{001}.
This means  that $C_1^2=\ol C_1^2 = -2$ on $S$, which implies that
$(C_1+ \ol C_1)\cdot C_1 = (C_1+\ol C_1)\cdot\ol C_1 = 0$.
Hence the restriction of the anticanonical bundle
 $\mathscr O_S(C_1+\ol C_1)$ on $S$ to the cycle $C_1 + \ol C_1$ is topologically trivial. 
From this it readily follows that
as a function of $l$, $h^0(lK_S^{-1})$ increases at most linearly as $l\to \infty$.
This implies that $h^0(Z, lF)$ increases at most quadratically as $l\to \infty$.
This contradicts our assumption that $Z$ is Moishezon.
Hence $m\neq 2$, as claimed.
\proofend

\vsp
Thus Moishezon twistor spaces on $4\CP^2$ with $h^0(F)=2$ can be classified 
by the number of  irreducible components of the base curve  of the pencil $|F|$, 
for which we can write by $2k$ where $2\le k\le 6$.

Let $\epsilon:S\to\CP^1\times\CP^1$ be the birational morphism as in the above proof.
Then by the real structure we can factorize $\epsilon$ as 
\begin{align}\label{dcp1}
\epsilon =\epsilon_4\circ\epsilon_3\circ\epsilon_2\circ\epsilon_1,
\end{align}
where each $\epsilon_i$ blows up a real pair of points on the anticanonical cycle.
Moreover, as we have denied the possibility $m=2$, even when we are in the situation \eqref{001},
some  $\epsilon_i$ has to blowup the pair of singular points of the curve
$C_1+\ol C_1$.
From this, by choosing  different blowdowns to $\CP^1\times\CP^1$,
 it follows that the situation \eqref{001} is absorbed in the situation \eqref{003}.
 Thus for any Moishezon twistor space on $4\CP^2$ with $h^0(F)=2$ 
we can  suppose that the image $\epsilon(C)$ is a $(2,2)$-curve consisting of 4 irreducible components.
Hence in the following we choose the birational morphism
$\epsilon$ so that the image $\epsilon(C)$ is 
a cycle of {\em four}\, rational curves.


We make use of the following easy property for $Z$, which is valid only on $4\CP^2$.
\begin{proposition}\label{prop:acs}
Let $Z$ be a Moishezon twistor space on $4\CP^2$ which satisfies $h^0(F)=2$.
Then if $S\in |F|$ is a smooth member,  we have 
\begin{align}\label{acs1}
h^0(2F) =  h^0(2K_S^{-1}) + 2.
\end{align}
\end{proposition}

\proof
As $H^1(F)=0$, from the exact sequence $0 \lra F \stackrel{s\otimes}{\lra} 2F \lra 2K_S^{-1}\lra 0$, where $s\in H^0(F)$  satisfies $(s)=S$ as before,
from $h^0(F)=2$ we get an exact sequence
\begin{align}\label{acs2}
0 \lra H^0(F) \stackrel{s\otimes}{\lra} H^0(2F) \lra H^0(2K_S^{-1}) \lra 0.
\end{align}
This implies \eqref{acs1}.
\proofend

\vsp

In the sequel we investigate structure of the anticanonical map of  $Z$ for each value of $k$,
by making use of the bi-anticanonical system of $S$ and the 
anticanonical system on $Z$.

\subsection{The case $k=6$}\label{ss:k6}
From the proof of Proposition \ref{prop:cycle},
$k$ attains the maximal
value $6$ iff any $\epsilon_i$ in the factorization \eqref{dcp1}  blows up a (real) pair of singular points of the cycle;
namely exactly when $S$ is a toric surface.
(Recall that the cycle $\epsilon(C)$ consists of 4 components.)

%
It is easy to see that we do not lose  any generality even if we suppose that  
$\epsilon_2\circ\epsilon_1$ blows up the 4 singular points of the cycle $\epsilon(C)$, 
and $\epsilon_3$ blows up one of the torus invariant 2 points on the strict transforms of $C_1$ and $\ol C_1$, where $C_1$ and $\ol C_1$ are curves in the expression \eqref{003}.
So the freedom for giving $S$ is solely in the choice of the pair of points blown-up by  $\epsilon_4$.
It is elementary to see that there are exactly 3 choices of $\epsilon_4$ which give
mutually non-isomorphic toric surfaces, and 
the sequence of  obtained by arranging the self-intersection numbers of irreducible components
of $C$ can be listed as follows:
\begin{align}
\label{string1}-4,   -1,   -2,   -2,   -2,   -1,   -4,   -1,   -2,   -2,   -2,   -1,\\
\label{string2}-3,   -2,   -1,   -3,   -2,   -1,  -3,   -2,   -1,  -3,   -2,   -1,\\
\label{string3}-3,    -1, -3,    -1,   -3,    -1,   -3,    -1,   -3,    -1,   -3,    -1.   
\end{align}
It is immediate to see that if the sequence for the cycle $C$ in $S$ 
is \eqref{string1}, then $h^0(K_S^{-1})=3$ which 
contradicts our assumption.
Hence it suffices to consider the other two types \eqref{string2} and \eqref{string2}.
These twistor spaces are investigated in \cite{Hon_JAG} under a strong assumption
that $Z$ has a $\mathbb C^*\times\mathbb C^*$-action.
We will show that the description obtained in \cite{Hon_JAG}
(in particular the explicit equations of the anticanonical model)
can be derived without assuming the existence of $\CC^*\times\CC^*$-action;
this makes the argument for obtaining projective models
of the twistor spaces of Joyce metrics \cite{J95} given in \cite[Section 2.2 and a part of Section 4]{Hon_JAG}
quite simpler:

%
%

\begin{proposition}\label{prop:k6}
Let $Z$ be a Moishezon twistor space on $4\CP^2$ satisfying $h^0(F)=2$.
Suppose that the number of irreducible components of 
the anticanonical cycle in $Z$ is $12$.
Then one of the following two situations  happens:
\begin{enumerate}
\item[(i)] 
If the self-intersection numbers for the cycle $C$ are as in \eqref{string2}, then $h^0(2F)= 7$ and  the anticanonical map $\Phi:Z\to \CP^6$ is birational to the image.
Further, the image is a complete intersection of 3 hyperquadrics, whose defining equations are of the form
\begin{align}\label{J1}
x_0x_1=x_2^2,\,\,
x_3x_4=q_1(x_0,x_1,x_2),\,\,
x_5x_6=q_2(x_0,x_1,x_2),
\end{align}
where $x_0,\cdots,x_6$ are homogeneous coordinates
on $\CP^6$ and $q_1$ and $q_2$ are quadratic 
homogeneous polynomials
in $x_0,x_1,x_2$.
\item[(ii)]
If the self-intersection numbers for the cycle $C$ are
as in \eqref{string3}, then
 $h^0(2F)=9$, and the anticanonical map $\Phi:Z\to \CP^8$ is birational to the image.
Further, the image is a 3-fold of degree $12$, which is not a complete intersection.
Furthermore, the defining equations of the image are
explicitly given as
\begin{align}\label{J2}
x_0x_1=x_2^2,\,\,
x_3x_4=q_1(x_0,x_1,x_2),\,\,
x_5x_6=q_2(x_0,x_1,x_2),\,\,
x_7x_8=q_3(x_0,x_1,x_2),\\
x_{2i+1}x_{2j+1} - x_0 x_{2k} = q_l(x_0,x_1,x_2),\,\,
x_{2i+2}x_{2j+2} - x_0 x_{2k+1} = q_m(x_0,x_1,x_2).\label{J2'}
\end{align}
where all $q_i$-s are quadratic polynomials
in $x_0,x_1,x_2$, and in \eqref{J2'} $i,j$ and $k$ take arbitrary combinations
satisfying $\{i,j,k\} = \{1,2,3\}$ (as sets).
\end{enumerate}
In particular, when $k=6$, the anticanonical map is
always birational.
\end{proposition}

\proof
We only give a proof for (i),
as (ii) can be proved (basically) in a similar way.
From \eqref{string2} we easily obtain $h^0(2K_S^{-1})=5$.
Hence by Proposition \ref{prop:acs} we get $h^0(2F)=7$.
It is also not difficult to see that the bi-anticanonical
map $\phi:S\to \CP^4$ is birational over the image.
Let $S^2H^0(F)$ denote the subspace of $H^0(2F)$ generated
by the image of the natural bilinear map $H^0(F)\times H^0(F) \to H^0(2F)$.
Then the rational map associated to the system 
$|S^2H^0(F)|$ factors as $Z\stackrel{\Phi_{|2F|}}{\lra} \mathbb PH^0(2F)^*
\stackrel{\pi}{\lra} \mathbb P\, S^2H^0(F)^*$, 
where $\pi$ is the projection associated to the inclusion
$S^2H^0(F)\subset H^0(2F)$,
and also as $Z\stackrel{\Phi_{|F|}}{\lra} \mathbb PH^0(2F)^*
\stackrel{g}{\lra} \mathbb P\, S^2H^0(F)^*$,
where $g$ is a natural map associated to the 
standard map $V^*\to S^nV^*$.  
Hence we have the  commutative diagram
\begin{equation}\label{016}
 \CD
 Z@>\Phi_{|2F|}>> \mathbb P H^0(2F)^*\\
 @V\Phi_{|F|} VV @VV\pi V\\
 \mathbb PH^0(F)^*@>g >> \mathbb P\, S^2H^0(F)^*.
 \endCD
 \end{equation}
This means that the anticanonical image $\Phi_{|2F|}(Z)$ is
always contained in $\pi^{-1}(g(\mathbb PH^0(F)^*)$.
If $h^0(F)=2$, $g$ is clearly an embedding 
of $\mathbb PH^0(F)^*=\CP^1$ as a conic.
Moreover, the diagram and the exact sequence \eqref{acs2}
also mean that the restriction of $\Phi_{|F|}$ 
to a general member $S\in |F|$ equals the bi-anticanonical
map $\phi$ of the surface.
From this the required birationality of $\Phi_{|2F|}$ 
follows.

For obtaining defining equations of the image  
$\Phi_{|2F|}(Z)$, 
by the method introduced in \cite[Section 2]{Hon-III},
we can explicitly find four members 
$Y_1,\ol Y_1,Y_2,\ol Y_2$ of the system 
$|2F|$ such that  any of their irreducible components are of degree one. 
(In \cite{Hon-III} the twistor spaces are assumed to have 
$\mathbb C^*\times\mathbb C^*$-action, but
for the purpose of finding these divisors
we do not need the action.
Also, these four members coincide with
the ones in \cite[Lemma 2.7]{Hon_JAG} which were
found rather incidentally.)

Let $x_3$ and $x_5$ be elements of $H^0(2F)$
such that $(x_3)=Y_1$ and $(x_5)=Y_2$,
and put $x_4=\ol {\sigma^*x_3}$
and $x_6=\ol {\sigma^*x_5}$.
Then $(x_4)=\ol Y_1$ and $(x_6)=\ol Y_2$.
Also, let $u_0, u_1$ be a basis of $H^0(F)$ and put 
$x_0=u_0^2,\,x_1= u_1^2$ and $x_2 = u_0u_1$.
Then $x_0,x_1,x_2$ form a basis of $S^2H^0(F)$ and 
the seven sections $x_0,\cdots,x_6$ form a basis of $H^0(2F)$.
Obviously we have the relation $x_0x_1=x_2^2$,
which is of course the equation of the image conic of $g$.
Next, since all irreducible components
of $Y_1$ are of degree one, $Y_1+ \ol Y_1$ 
is a sum of four members of the pencil $|F|$,
meaning $x_3x_4\in S^4H^0(F)$.
By the same reason (or the real structure) we have $x_5x_6\in S^4H^0(F)$.
Therefore we can write
\begin{align}\label{key1}
x_3x_4 = a_1u_0^4+ a_2 u_0^3 u_1 + a_3 u_0^2 u_1^2
+ a_4 u_0 u_1^3 + a_5 u_1 ^4.
\end{align}
Then the right-hand-side of \eqref{key1}
can be rewritten as a quadratic polynomial
of $x_0,x_1,x_2$.
(The last quadratic polynomial is never unique, 
as we have the relation $x_0x_1=x_2^2$.)
This gives the second equation in \eqref{J1}.
The same argument clearly works for $x_5x_6$,
which gives the third equation in \eqref{J1}.
Then it is not difficult to see that the threefold
defined by the three equations \eqref{J1}
is irreducible.
Hence the anticanonical image $\Phi_{|2F|}(Z)$ is
defined by the equations \eqref{J1}.
\proofend

\vsp
By comparing the equations \eqref{J1}--\eqref{J2'} with those  obtained in
 \cite[Theorems 2.3 and  4.4]{Hon_JAG}, 
the defining equations of the anticanonical models
 coincide with the case for the twistor spaces
of Joyce metrics.
So it seems very natural to expect
 that if a twistor space on $4\CP^2$
(or $n\CP^2$, more generally)
has an irreducible real member $S\in |F|$ which is a toric surface, then $Z$ is a twistor space of Joyce metric.
(Of course the point here is that we are not assuming
the existence of torus-action.)
This might give a new characterization of Joyce metrics.

\subsection{The case $k=5$}
These twistor spaces can be classified into 2 separate types as follows:
\begin{proposition}\label{prop:k5}
Let $Z$ be a Moishezon twistor space on $4\CP^2$ satisfying $h^0(F)=2$.
Suppose the number of irreducible components
of the anticanonical cycle $C$ in $Z$ is $10$.
Then one of the following two situations  happens:
\begin{enumerate}
\item[(i)] $h^0(2F)= 7$ and  the anticanonical map $
\Phi:Z\to \CP^6$ is birational to the image;
further, the image is a complete intersection of 3 hyperquadrics.
\item[(ii)] $h^0(2F)=5$, and the image of the anticanonical map $\Phi:Z\to \CP^4$ is  a scroll of planes over a conic, and
the anticanonical map is 2 to 1 over the scroll.
Moreover restriction of the branch divisor of $\Phi$ to a general 2-plane
of the scroll is a quartic curve.
\end{enumerate}
\end{proposition}

\proof Let $\epsilon:S\to\CP^1\times\CP^1$ be the birational morphism as in Section \ref{ss:bc} and again decompose  as $\epsilon =\epsilon_4\circ\epsilon_3\circ\epsilon_2\circ\epsilon_1$,
where each $\epsilon_i$ blows up a real pair of points belonging to the anticanonical cycle.
Further, as is already mentioned, we can suppose that the image $\epsilon(C)$ on
$\CP^1\times\CP^1$ is a cycle of 4 curves (as in \eqref{003}). 
Then from the assumption $k=5$, 
 there exists a unique $j$ such that $\epsilon_j$ blows up a pair of {\em smooth} points  on the anticanonical cycle.
We rearrange the order of $\epsilon_i$-s in a way that $j=4$ holds.
Put $S_1:=\epsilon_4(S)$.
By construction $S_1$ has an anticanonical cycle which consists of 10 irreducible components, and it is a toric surface.
As mentioned in Section \ref{ss:k6}, the possible structure of $S_1$ is unique,
and the sequence of the self-intersection numbers
of  irreducible components of the anticanonical cycle
(on $S_1$) in order is
\begin{align}\label{3563}
-3,-1,-2,-2,-1,-3,-1,-2,-2,-1.
\end{align}
Let $C_1,\cdots,C_5,\ol C_1,\cdots,\ol C_5$ be the  components arranged in this order (so  $C_1^2=\ol C_1^2=-3$ for instance.)
Note that $S_1$ has an involution which preserves the real structure and 
which exchanges $C_1,C_2,C_3,C_4$ and $ C_5$
with $\ol C_1, C_5,C_4,C_3$ and $C_2$ respectively.
Hence we can suppose that $S$ is obtained from $S_1$ by blowing-up a smooth point on $C_1,C_2$ or $C_3$ and its conjugate point,
where a smooth point refers a point which is not a
singularity of the cycle.
But if  $\epsilon_4$ blows up a point on $C_1$, then $S$ satisfies $h^0(K_S^{-1})=3$ and this contradicts  $h^0(K_S^{-1})=1$.
So $\epsilon_4$ blows up a  point of $C_2$ or a  point on $C_3$, which is not on the 
singularities of the anticanonical cycle.

Suppose that $\epsilon_4$ blows up a smooth point of $C_2$.
Then the structure of $S$ is the same as that of fundamental divisors
in the twistor spaces studied in \cite{HonDSn1} (in the case of $4\CP^2$, of course),
 and hence by \cite[Proposition 2.1 (ii) and (v)]{HonDSn1}
the bi-anticanonical system of $S$ is 2-dimensional, whose associated map is a degree 2 morphism onto $\CP^2$.
In particular by Proposition \ref{prop:acs} we obtain $h^0(2F) = 5$.
Further, since the restriction of $\Phi_{|2F|}$ to any smooth member $S\in |F|$ 
is exactly the bi-anticanonical map by the exact sequence \eqref{acs2},
from the diagram \eqref{016}
we conclude that $\Phi_{|2F|}$ is 2 to 1 over the scroll of planes over a conic.
Furthermore, we readily see that the branch divisor of the bi-anticanonical map
of $S$ is a quartic curve.
Thus $Z$  satisfies all the properties of the case (ii) in the proposition.

Next we suppose that $\epsilon_4$ blows up a  point of $C_3$.
In this case, from \eqref{3563}, the sequence of self-intersection numbers of the anticanonical cycle 
on $S$ becomes
\begin{align}
-3,-1,-3,-2,-1,-3,-1,-3,-2,-1.
\end{align}
For this surface $S$ we have the following

\begin{lemma}\label{lemma:83}
The bi-anticanonical system of $S$ satisfies the following.
\begin{enumerate}
\item
[(i)] $h^0(2K_S^{-1})=5$.
\item
[(ii)]  The fixed components of the system $|2K_S^{-1}|$ is $C_1 + C_3 + C_4 
+ \ol C_1 + \ol C_3 + \ol C_4$, and after removing this,
the system becomes base point free.
\item
[(iii)] If $\phi:S\to \CP^4$ is the
morphism induced by the bi-anticanonical system, then $\phi$ is birational over the image  $\phi(S)$, and 
$\phi(S)$
is a complete intersection of two hyperquadrics in $\CP^4$.
\end{enumerate}
\end{lemma}

\proof
The curves $C_1, C_3, C_4$ and their conjugations can be readily seen to be 
fixed components from computation of the intersection numbers.
Let $B$ be the sum of these 6 curves (with each coefficient being one).
Then we easily have $(2K_S^{-1}-B)\cdot C_i=(2K_S^{-1} - B)\cdot \ol C_i= 0$
for $i=2,5$.
Hence noticing $(2K_S^{-1}-B) -( C_2 +C_5 + \ol C_2 + \ol C_5 ) \simeq K^{-1}_S$
 we obtain an exact sequence
\begin{align}\label{ses:25}
0 \lra K_S^{-1} \lra 2K_S^{-1}
\otimes\mathscr O_S(-B) \lra 
\mathscr O_{C_2 \,\cup\, C_5 \,\cup\, \ol C_2 \,\cup\, \ol C_5} \lra 0,
\end{align}
where the second arrow is a multiplication of a defining section of the divisor
$C_2 +C_5 + \ol C_2 + \ol C_5$.
Further we readily have $h^0(K_S^{-1})=1$ and hence $h^1(K_S^{-1})=0$ by Riemann-Roch.
Therefore the cohomology exact sequence of \eqref{ses:25} gives $h^0(2K_S^{-1}) = 5$,
and also implies that ${\rm{Bs}}\, |2K_S^{-1} - B|$ is disjoint from $C_2 \cup C_5 \cup \ol C_2 \cup \ol C_5$. 
Hence, as $(2K_S^{-1} - B) \cdot C_3 =0$, 
 ${\rm{Bs}}\, |2K_S^{-1} - B|$ is disjoint from $C_3\cup \ol C_3$ also.
 Therefore  if ${\rm{Bs}}\, |2K_S^{-1} - B|$ is non-empty, it must be points on 
 $C_1\cup C_4 \cup \ol C_1 \cup \ol C_4$ which are not singularities of the anticanonical cycle on $S$.
 But this cannot happen since $S$ has a non-trivial $\mathbb C^*$-action which fixes any
 point on  $C_3\cup \ol C_3$, and this $\CC^*$-action automatically lifts on the system $ |2K_S^{-1} - B|$,
 so that its base locus must be $\CC^*$-invariant.
 Hence we obtain (i) and (ii).
 
 To show (iii), we notice that the two divisors
 $$
 f:= C_4 + 2 C_5 + \ol C_1 + \ol C_2 
 \quad{\text{and}}\quad
\ol f:=\ol  C_4 + 2 \ol  C_5 +  C_1 +  C_2
 $$
 are mutually linearly equivalent, and that $2C-f+\ol f$ and $2C + f - \ol f$ are effective divisors,
 where $C$ is the unique anticanonical curve (cycle) on $S$ as before. 
 Therefore it makes sense to speak about the linear subsystem of $|2K_S^{-1}|$ generated by 
 \begin{align}
 \label{3curves}
 2C-f+\ol f,\,\, 2C + f - \ol f
  \quad{\text{and}}\quad
  2C.
 \end{align}
Since $B=\Bs\,|2K_S^{-1}|$, the curve $B$
 contained in  the fixed components of this subsystem.
 By looking at the multiplicities of the components, we readily see that 
 these 3 divisors are linearly independent, so that 
 they determine a 2-dimensional subsystem of $|2K_S^{-1}|$.
 Write $V$ for the 3-dimensional linear subspace of $H^0(2K_S^{-1})$
 corresponding to this 2-dimensional subsystem.
 Then in a similar way to the diagram \eqref{016} we have the commutative diagram
\begin{equation}\label{019}
 \CD
 S@>\phi>> \CP^4\\
 @VVV @VV{p}V\\
 \CP^1@>{\iota}>> \CP^2
 \endCD
 \end{equation}
where $\phi$ is the morphism associated to 
$|2K_S^{-1}|\simeq |2K_S^{-1} - B|$, and the composition $S\to \CP^1\to \CP^2$ is the rational map associated to the 2-dimensional subsystem $|V|$,
and $p$ is a projection associated to the subspace $V\subset H^0(2K_S^{-1})$.
Since the sum of the first two curves in \eqref{3curves} equals the twice of the last curve,
the image of the   rational  map associated to $|V|$ is a conic.
This diagram implies that the bi-anticanonical image $\phi(S)$ is 
contained in   $p^{-1}(\iota(\CP^1))$, the scroll of 2-planes over a conic.
Hence in order to show that $\phi(S)$ is a complete intersection of two hyperquadrics,
it is enough to show that the image under $\phi$ of  general fibers of $S\to \CP^1$
is a conic.
For this we note that the movable part of the 2-dimensional linear subsystem $|V|$ is composed with the pencil $\langle f, \ol f \rangle$.
Let $D$ be a general member of this pencil, which is easily seen to be a smooth rational 
curve.
By computing self-intersection numbers, we readily deduce that 
the system $|2K_S^{-1}-B-D|$ has $|f|$ 
as a movable part.
Further the last linear system is a pencil.
Hence from the exact sequence
$$
0 \lra 2K_S^{-1}\otimes\mathscr O_S(-B-D) \lra 
2K_S^{-1}\otimes\mathscr O_S(-B) \lra \mathscr O_D(2) \lra 0
$$
we obtain that the image of the restriction map
$H^0(2K_S^{-1}\otimes\mathscr O_S(-B))\to H^0(\mathscr O_D(2))$ is surjective.
This means that the image $\phi(D)$ is a conic.
Hence $\phi(S)$ is an intersection of two hyperquadrics.
This also implies
that the restriction $\phi|_D:D\to \phi(D)$ is biholomorphic.
Then by the diagram \eqref{019},  $\phi$ is birational over the image $\phi(S)$,
finishing a proof of the lemma.
\proofend

\begin{remark}
{\em
The bi-anticanonical map $\phi:S\to\CP^4$ contracts $C_2, C_3$ and $C_5$,
and the image $\phi(S)$ has singularities.
}
\end{remark}

\noindent
{\em Continuation of the proof of Proposition \ref{prop:k5}.}
By Lemma \ref{lemma:83} (i) and Proposition \ref{prop:acs} we have $h^0(2F) = 7$.
Further the diagram \eqref{016} is valid also in the present situation,
and the dimensions of the projective spaces in
the diagram are precisely the same as the case (ii)
of the proposition.
Also, the image of $g\circ\Phi_{|F|}$ is still a conic in $\mathbb P S^2H^0(F)^*=\CP^2$.
Hence by the same reason to the case (ii), the image $\Phi_{|2F|}(Z)$ is
contained in the scroll of planes over the conic.
Moreover, again by the surjectivity of the restriction map in the short exact sequence
\eqref{acs2}, the restriction of $\Phi_{|2F|}$ to a smooth member $S$ of $|F|$ is
exactly the bi-anticanonical map of $S$.
Hence in order to prove that the image $\Phi_{|2F|}(Z)$ is a complete intersection of three hyperquadrics
(in $\CP^6$), it suffices to show that a general fiber of $\Phi_{|F|}$ is mapped onto
a complete intersection of two hyperquadrics in a fiber of $\pi$.
This is exactly the assertion (iii) of Lemma \ref{lemma:83}.
Finally, the birationality of $\Phi_{|2F|}$ now follows immediately from the diagram 
\eqref{016} and the birationality of $\phi$
proved as Lemma \ref{lemma:83} (iii).
Thus we are  exactly in the situation (i) of Proposition \ref{prop:k5}.
\proofend

\vsp
We note that the twistor spaces studied in \cite{HonDSn1} (in the case of $4\CP^2$)
fall into the case (ii) in
Proposition \ref{prop:k5}.
On the other hand, the twistor spaces in the case (i) seem to have 
not appeared in the literature.
But it is quite likely that these twistor spaces can be 
considered as mild degenerations of the twistor spaces
constructed in \cite{Hon-I}.

\subsection{The case $k=4$}
These twistor spaces can be classified into three  types as follows:
\begin{proposition}\label{prop:k4}
Let $Z$ be a Moishezon twistor space on $4\CP^2$ satisfying $h^0(F)=2$.
Suppose the number of irreducible components of the
 anticanonical cycle of $Z$ is $8$.
Then one of the following three situations  happens:
\begin{enumerate}
\item[(i)] $h^0(2F)=7$,  the anticanonical map
$\Phi:Z\to\CP^6$ is birational, and 
the image $\Phi(Z)$ is a complete intersection of three hyperquadrics.
\item[(ii)] $h^0(2F)=5$, and the image of the anticanonical map $\Phi:Z\to \CP^4$ is  a scroll of planes over a conic, and
$\Phi$ is rationally 2 to 1 over the scroll.
Moreover restriction of the branch divisor of $\Phi$ to a general 2-plane
of the scroll is a quartic curve.
\item[(iii)] $h^0(2F)= 5$ and the image of the anticanonical map $Z\to \CP^4$ is a complete
intersection of 2 hyperquadrics, whose defining equations
are of the form
\begin{align}\label{ci}
x_0x_1 = x_2^2,\,\,x_3x_4 = q(x_0,x_1,x_2),
\end{align} 
where $x_0,\cdots,x_4$ are homogeneous coordinates on
$\CP^4$ and $q$ is a quadratic polynomial in $x_0,x_1,x_2$.
\end{enumerate}
\end{proposition}

\begin{remark}
{\em As will be clear at the end of Section \ref{ss:k32},
the case (iii) is the unique situation where
the  image of the anticanonical system is not 3-dimensional,
provided that $h^0(F)=2$.
In other words, except this case, the anticanonical image is 3-dimensional when $h^0(F)=2$.
}
\end{remark}

\noindent
{\em Proof of Proposition \ref{prop:k4}.}
As in the cases of $k=6$ and $k=5$, take a real member $S\in|F|$ and take the 
birational morphism $\epsilon:S\to \CP^1\times\CP^1$ such that 
$\epsilon(C)$ consists of 4 components (as in \eqref{001}).
Again we factorize 
$\epsilon$ as $\epsilon =\epsilon_4\circ\epsilon_3\circ\epsilon_2\circ\epsilon_1$. 
 Then by the assumption 
$k=4$ we can suppose that $\epsilon_2\circ\epsilon_1$  blows up the 4 singular points of the 
$(2,2)$-curve $\epsilon(C)$, and that $\epsilon_3$ and $\epsilon_4$ blowup smooth points 
of the anticanonical cycle.
Let $S_1$ be the surface obtained from $\CP^1\times\CP^1$ by applying $\epsilon_2\circ\epsilon_1$,
and we name the components of the anticanonical cycle on $S_1$ as
$C_1,C_2,C_3,C_4,\ol C_1,\ol C_2,\ol C_3,\ol C_4$ in a way that adjacent curves intersect as before and $C_1$ is a $(-2)$-curve.
By the real structure one of the 2 blown-up points of $\epsilon_3$ has to be
on $C_1\cup C_2\cup C_3\cup C_4$, and the same for  $\epsilon_4$.
Hence $S$ can be specified by choosing 2 points on $C_1\cup C_2\cup C_3\cup C_4$, which are not  the singularities
 of the cycle.
In the following we do not mention the conjugate 
blown-up points.

First suppose that $\epsilon_4\circ\epsilon_3$ blows up 2 points on $C_1$;
here we are allowing these points to be a same point, in which case $\epsilon_4$ blows up
the intersection point of the exceptional curve with the strict transform of $C_1$.
Including this case it is readily seen that $h^0(K_S^{-1})=3$, contradicting $h^0(F)=2$.
Hence by an obvious symmetry of $S_1$, we also
obtain that  $\epsilon_4\circ\epsilon_3$ cannot
blowup 2 points on $C_3$.

Next suppose that $\epsilon_4\circ\epsilon_3$ blows up one point on $C_1$ and one point on $C_2$.
Then using the same letter to denote the strict transforms, the self-intersection numbers of the 8 curves $C_1,C_2,\cdots,\ol C_4$ on $S$ respectively become
\begin{align}
-3,-2,-2,-1,-3,-2,-2,-1.
\end{align}
From this by standard computations we can deduce that the system $|2K_S^{-1}|$ is 2-dimensional with the fixed components 
$C-C_4-\ol C_4$, and after removing it the system becomes base point free,
which induces a degree 2 morphism $S\to\CP^2$ whose branch
divisor is a quartic curve.
Hence again by
Proposition \ref{prop:acs} and the diagram \eqref{016}, we obtain that $Z$ satisfies the property  (ii) of the proposition.
By a symmetry on $S_1$ the same conclusion holds for the case that 
$\epsilon_4\circ\epsilon_3$ blows up one point on $C_1$ and one point on $C_4$,
and also the case that it blows up one point on $C_2$ and one point on $C_3$.

Next suppose that $\epsilon_4\circ\epsilon_3$ blows up one point on $C_1$ and one point on $C_3$.
Then this time the self-intersection numbers of the 8 components $C_1,C_2,\cdots,\ol C_4$ 
of the anticanonical cycle $C$ on $S$ respectively become
\begin{align}\label{3131}
-3,-1,-3,-1,-3,-1,-3,-1.
\end{align}
Note that unlike the case $k=5$, this surface $S$ does not have a
non-trivial $\mathbb C^*$-action.
\begin{lemma}\label{lemma:bas2}
The bi-anticanonical system of this surface $S$ has the following properties:
\begin{enumerate}
\item
[(i)] $h^0(2K_S^{-1})=5$, 
\item
[(ii)] the fixed components of \,$|2K_S^{-1}|$ is $C_1+C_3+\ol C_1+\ol C_3$,
\item
[(iii)] after removing the fixed components, the system becomes base point free,
\item
[
(iv)] the morphism $\phi:S\to\CP^4$ induced by the system contracts 
the four $(-1)$-curves $C_2, C_4, \ol C_2, \ol C_4$,
and the image $\phi(S)$ is a quartic surface.
\item
[(v)] if $\phi$ does not contract any other curves, then a general hyperplane section of the quartic surface $\phi(S)$ 
is a smooth elliptic curve.
\end{enumerate}
\end{lemma}

Note that in (iv) we are not claiming that the birational morphism
$\phi$ does not contract any curve other than the four
$(-1)$-curves.
Actually if there exists a $(1,1)$-curve
on $\CP^1\times\CP^1$ going through the 4 
blown-up points of $\epsilon_4\circ\epsilon_3$,
then its strict transform is contracted to 
an ordinary double point.
Later we will show that 
it cannot happen that this kind of curve exists
for {\em any} smooth member $S$ of the pencil $|F|$.
\proof
For (i)--(iii), as $K_S^{-1} \cdot C_i=K_S^{-1}\cdot \ol C_i=-1$ for $i=1,3$, these curves are
base curves of $|2K_S^{-1}|$.
Put $B:= C_1+C_3+\ol C_1+\ol C_3$, and let $s\in H^0(\mathscr O_S(C_2+C_4+\ol C_2 + \ol C_4))$ 
be an element  satisfying $(s) = C_2+C_4+\ol C_2 + \ol C_4$.
Then by the exact sequence 
\begin{align}\label{ses:4}
0 \lra K_S^{-1} \stackrel{\otimes s}{\lra} 2K_S^{-1} \otimes\mathscr O_S(-B) \lra \mathscr O_{C_2 \,\cup\, C_4\,\cup\,\ol C_2 \,\cup\, \ol C_4}\lra 0
\end{align}
and noticing $H^1(K_S^{-1})=0$, we obtain $h^0(2K_S^{-1})=5$ and that $|2K_S^{-1}
-B|$ has no fixed component, meaning (i) and (ii).
Suppose ${\rm Bs}\,|2K_S^{-1}-B|\neq\emptyset$.
Then again by the cohomology exact sequence of \eqref{ses:4}, 
it is contained in the support of $B$, and also it is disjoint from $C_2\cup C_4\cup\ol C_2 \cup \ol C_4$.
Noting that the restriction of the line bundle 
$2K_S^{-1} \otimes\mathscr O_S(-B)$ to the $(-3)$-curves $C_i$ and $\ol C_i$ 
is of degree 1, base points appears on each of these curves at one point at most.
By reality, the number of base points is either 2 or 4.
If it is 2, the base points are resolved by blowing-up the 2 points,
and consequently we obtain a morphism $\tilde S\to \CP^4$, where
$\tilde S$ is the 2 points blowup.
Then since $(2K_S^{-1}-B)^2 = 4$ on $S$, the self-intersection number of the 
base point free linear system on $\tilde S$ is 2.
This means that the image of $\tilde S$ in $\CP^4$ is quadratic surface.
But a quadratic surface in $\CP^4$ necessarily degenerates (see \cite{R}),
and hence the number of the base points (of $|2K_S^{-1}-B|$) cannot be 2.
If the number of the base points is 4,
 $|2K_S^{-1}-B|$ has one base point on each of the four $(-3)$-curves,
and by blowing-up these 4 points we obtain a 
surface $\tilde S$ and a fixed point free linear system on $\tilde S$. 
Since the self-intersection number of this system is zero, the image of the morphism $\tilde S\to \CP^4$ 
induced by the system is a curve.
Let $E$ be any one of the 4 exceptional curves of $\tilde S\to S$.
Then the intersection number of $E$ with the free system on $\tilde S$ is 1.
Hence the image of $E$ under the morphism $\tilde S\to \CP^4$ must be a line.
These mean that the image curve in $\CP^4$ is a line, which cannot happen.
Thus the number of the base points of $|2K_S^{-1}-B|$ cannot be 4.
This implies that  $|2K_S^{-1}-B|$ is base point free.
Hence we obtain (iii).

Let $\phi:S\to\CP^4$ be the morphism induced by $|2K_S^{-1}-B|$.
Since $(2K_S^{-1}-B)^2=4$, the image $\phi(S)$ is a quadratic surface and $\phi$ is of degree 2, or
otherwise $\phi(S)$ is a quartic surface and $\phi$ is birational.
But again it cannot be quadratic, so $\phi(S)$ is birational and $\phi(S)$ is a quartic surface.
Then as $(2K_S^{-1}-B)\cdot C_i = (2K_S^{-1}-B)\cdot\ol C_i = 0$ for $i=2,4$,
these four $(-1)$-curves are contracted to points by $\phi$.
Thus we get (iv) of the lemma.
Finally, if $\phi$ does not contract any curve except the four $(-1)$-curves, 
$\phi(S)$ is non-singular.
Further, by Bertini's theorem, a general hyperplane section of $\phi(S)$
is irreducible and non-singular. 
As such a hyperplane section has to be a biholomorphic image of
a member of $|2K_S^{-1}-B|$, in order to prove (v), it suffices to compute the
arithmetic genus of the system, which is given by
\begin{align*}
\frac12 (2K_S^{-1}-B) (2K_S^{-1}-B+K_S)+1 & =
\frac12 (2K_S^{-1}-B)^2 + \frac12 (2K_S^{-1}-B)K_S+1\\
&= \frac12 4 + 0 + \frac12 (-4) + 1 = 1.
\end{align*}
Hence we obtain (v) of the lemma.
\proofend

\vsp
For finishing a proof of Proposition \ref{prop:k4} 
we  need another lemma about the anticanonical map 
of the twistor space $Z$.
Recall that we are considering a twistor space $Z$
such
that $\epsilon_4\circ\epsilon_3$ (for a real irreducible member $S\in |F|$ on $Z$) blows up
one point on $C_1$ and one point on $C_3$.
From the fact that $C_1$ and $C_3$ are $(-3)$-curves on $S$, 
it is easy to see that the normal bundle satisfies $N_{C_i/Z}\simeq
\mathscr O(-2)^{\oplus 2}$ for $i=1,3$, and hence 
if we blowup these curves, each exceptional divisor is 
biholomorphic to the product $\CP^1\times\CP^1$.
The following lemma clarifies the structure of the anticanonical map
of the twistor space $Z$
 to a considerable degree:
\begin{lemma}\label{lemma:nc}
Let $\mu:\tilde Z\to Z$ be the blow-up at the base curve $B=C_1\cup C_3\cup \ol C_1\cup\ol C_4$, $E_i$ and $\ol E_i$ $(i=1,3)$  the exceptional divisors, and put $E=E_1+E_3+\ol E_1+\ol E_3$.
(The 4 curves are of course disjoint.)
Then the system $|\mu^*(2F)-E|$ is base point free and 
the resulting morphism contracts $E_i\simeq\CP^1\times\CP^1$
 and $\ol E_i$ $(i=1,3)$
to $\CP^1$ along the projection which is different from the original projections
$E_i\to C_i$ and $\ol E_i\to \ol C_i$.
Moreover, other than these 4 exceptional divisors, the morphism does not contract any divisors (in $\tilde Z$) to a lower-dimensional subvariety.  
\end{lemma}

\proof
By the exact sequence \eqref{acs2}  
we have $\Bs\, |2F| = \Bs\, |2K_S^{-1}|$ and 
the latter is the curve $B$ by Lemma \ref{lemma:bas2} (ii).
Moreover by the exact sequence \eqref{acs2},
a general member of $|2F|$ intersects $S$ transversally along $B$.
This means $\Bs\,|\mu^*(2F)-E|=\emptyset$. 
For any of the 4 exceptional divisors, let $\mathscr O(0,1)$ be
a fiber-class of the projection $\mu$.
Then we have $N_{E_i/\tilde Z}\simeq\mathscr O(-1,-2)$.
Hence we obtain
\begin{align}\label{nb2}
(\mu^*(2F)-E)|_{E_i}
\simeq
\mu^*(2K_S^{-1}|_{C_i})\otimes N_{E_i/\tilde Z}^{-1}
\simeq
\mathscr O(0,-2) \otimes \mathscr O(1,2)
\simeq \mathscr O(1,0). 
\end{align}
If $\tilde{\Phi}:\tilde Z\to \CP^6$ denotes the morphism
associated to $|\mu^*(2F)-E|$, this implies that 
$\tilde{\Phi}$ contracts $E_i$ to $\CP^1$ along the projection
which is different from the original $E_i\to C_i$.

It remains to see that the morphism $\tilde{\Phi}$ does not
contract any divisor other than the components of $E$.
This can be proved in a similar way to \cite[Proposition 3.6]{HonDS4_1},
so here we only present an outline.
It is enough to see that there does not exist a real
irreducible divisor $D$ on $\tilde Z$ with $D\neq E_i,\ol E_i$
such that $(\mu^*(2F)-E)^2\cdot D=0$.
Such a divisor is linearly equivalent to $\mu^*(kF) - l_1(E_1+ \ol E_1) - 
l_3 (E_3 + \ol E_3)$ for some $k\ge 1$ and $l_1\ge 0, l_3\ge 0$.
But a computation shows that the intersection number of 
the last class with $(\mu^*(2F) - E)^2$ is $4k$.
Hence such a divisor $D$ does not exist.
\proofend

\vsp
We note that from the normal bundle of $E_i$ in $\tilde Z$ and also from
\eqref{nb2}
the image $\tilde{\Phi}(E_i)$ and $\tilde{\Phi}(\ol E_i)$ ($i=1,3$) are lines in $\CP^6$ and the anticanonical image $\Phi(Z)=\tilde{\Phi}(Z)$ has ordinary double points along these 4 lines.
We also note that $\tilde{\Phi}$ contracts not only the exceptional divisors of $\mu$ but also the strict transforms of the remaining
curves $C_i$ and $\ol C_i$ for $i=2,4$.
With these informations, it might be interesting to find explicit form of quadratic polynomials
which defines the anticanonical model.
(We do not know whether there is a curve 
which are contracted by $\tilde{\Phi}$ except the above four curves.)

\vsp
\noindent
{\em Continuation of the proof of Proposition \ref{prop:k4}.} 
Recall the situation that $Z$ has $S\in |F|$ whose unique anticanonical curve
is presented as in \eqref{3131}.
We claim that for a general choice of a real member $S\in|F|$, the bi-anticanonical
map $\phi:S\to \CP^4$ does not contract any curve other than
the four $(-1)$-curves $C_2,C_4,\ol C_2$ and $\ol C_4$.
If this is not the case, then from the diagram \eqref{016} and  since the restriction of the 
anticanonical map $\Phi:Z\to \CP^6$ to $S$ is exactly
the bi-anticanonical map,
$\tilde{\Phi}$ must contract a {\em divisor}, say $D$,  to a lower-dimensional subvariety.
But by Lemma \ref{lemma:nc} such a divisor $D$ 
must be one of $E_i$ and $\ol E_i$ $(i=1,3)$.
Thus we conclude that for a general member $S\in |F|$
the bi-anticanonical map $\phi$ is just
a (simultaneous) blowdown of the $(-1)$-curves
$C_2,C_4,\ol C_2$ and $\ol C_4$.
This means that for these surfaces the image $\phi(S)\subset\CP^4$ 
is a {\em non-singular} surface.
The latter surface is of degree 4 by Lemma \ref{lemma:bas2} (iv).
It is classically known (see \cite{R}) that 
 a smooth quartic surface in $\CP^4$ is either
the image of 
a Veronese surface under a projection  $\CP^5\to \CP^4$,
or a (complete) intersection of two hyperquadrics,
and that these two quartic surfaces are distinguished by
the genus of a generic hyperplane section;
for the Veronese surface the genus is zero,
while  for the complete intersection of the quadrics
the genus is one.
Hence by Lemma \ref{lemma:bas2} (v) we can conclude that 
$\phi(S)$ is an intersection of two hyperquadrics.
Thus we obtain that for a generic
member of the pencil $|F|$, 
 $\phi(S)$ is a complete intersection of two hyperquadrics.
Hence so is $\Phi(S)$.
Then by the diagram \eqref{016}, we obtain that 
$\Phi(Z)$ is a complete intersection of three
hyperquadrics. Thus we have seen
that $Z$ satisfies all the  properties 
stated in (i) of  Proposition \ref{prop:k4}.

Next suppose that $\epsilon_4\circ\epsilon_3$ blows up two points on $C_2$.
Then this time the sequence of the self-intersection numbers of the 8 curves $C_1,C_2,\cdots,\ol C_4$ on $S$ become
\begin{align}
-2,-3,-2,-1,-2,-3,-2,-1.
\end{align}
It follows that the system $|2K_S^{-1}|$ is 2-dimensional with the fixed components
$C_1+2C_2+C_3+\ol C_1+ 2\ol C_2+\ol C_3$, and that after removing it the system becomes base point free.
Moreover, the last system is generated by
 a pencil $|C_3 + 2 C_4 + \ol C_1| = 
|\ol C_3 + 2 \ol C_4 +  C_1|$, and hence
it  gives a morphism to a conic.
Therefore the diagram \eqref{016} means that the image $\Phi_{|2F|}(Z)$ is 2-dimensional, and 
it is a complete intersection of two hyperquadrics.
So we are in the situation (iii).
The claim about the equations \eqref{ci} of the two hyperquadrics
can be derived by the same way to the case $k=6$ given 
in the proof of Proposition \ref{prop:k6} and 
we omit the detail.
For the claim about fibers of $\Phi$, it suffices to show that general fibers 
of the bi-anticanonical map  $\phi$ of 
$S$ are smooth rational curves.
But as above the movable part of $|2K_S^{-1}|$ is
generated by the pencil 
$|C_3 + 2C_4 + \ol C_1|$,
and it is easy to see that general members of
this pencil are  non-singular rational curves.
Hence general fibers of $\phi$ are non-singular rational curves.
Thus we have obtained all the properties in (iii) of the proposition.

If $\epsilon_4\circ\epsilon_3$ blows up one point on $C_2$ and one point on $C_4$,
all the components of $C$  become $(-2)$-curves.
In this case the line bundle $K^{-1}|_C$ clearly becomes topologically trivial, and it follows that 
$h^0(lK_S^{-1})$ increases at most linearly about $l$.
This means that $Z$ is non-Moishezon.
By symmetries of $S_1$ we do not need to consider the remaining cases, because they are reduced to
one of the  above cases.
Thus we have completed a proof of Proposition \ref{prop:k4}.
\proofend

\vsp
We note that the twistor spaces appeared  in 
\cite[Proposition 5.4, Figure 5 (iv)]{HonDS4_1}
belong to the situation (i) of Proposition \ref{prop:k4}, and the twistor spaces studied in \cite{Hon-I} belong to the case (iii) of Proposition \ref{prop:k4}.
On the other hand, the twistor spaces in the case (ii) seem to have not appeared in the literature.

\subsection{The cases  $k=3$ and $k=2$}\label{ss:k32}
The conclusion in these two cases is  simple:
\begin{proposition}\label{prop:k32}
Let $Z$ be a Moishezon twistor space on $4\CP^2$ satisfying $h^0(F)=2$.
Suppose the number of irreducible components of the
 anticanonical cycle in $Z$ is $6$ or $4$.
Then $h^0(2F)=5$, and the image of the anticanonical map $Z\to \CP^4$ is  a scroll of planes over a conic, and
the anticanonical map is 2 to 1 over the scroll.
Moreover restriction of the branch divisor of $\Phi$ to a general 2-plane
of the scroll is a quartic curve.
\end{proposition}

\proof
Let $S\in |F|$ and  $\epsilon =\epsilon_4\circ\epsilon_3\circ\epsilon_2\circ\epsilon_1
:S\to\CP^1\times\CP^1$
be as in the proof of Propositions \ref{prop:k6},
\ref{prop:k5} or \ref{prop:k4}.
If $k=3$, we can suppose that $\epsilon_1$ blows up a pair of singularities of the cycle $\epsilon(C)$ and 
the remaining $\epsilon_i$-s blow up  pairs of smooth points of the cycle.
Like before let $\epsilon_1$ be the surface obtained
by applying $\epsilon_1$, and
$C_1,\cdots,\ol C_3$  the cycle of six $(-1)$-curves on $S_1$.
If $\epsilon_4\circ\epsilon_3\circ\epsilon_2$ blows up three points on $C_1$,
we get $h^0(K_S^{-1})=3$, which contradicts our assumption.
If $\epsilon_4\circ\epsilon_3\circ\epsilon_2$ blows up two points on $C_1$,
by a symmetry of $S_1$ we can suppose that it blows up one point on $C_2$, and the self-intersection numbers of the 
anticanonical cycle  $C$ on $S$ become
\begin{align}
-3,-2,-1,-3,-2,-1.
\end{align}
From this we deduce that the system $|2K_S^{-1}|$ is 2-dimensional, the fixed components
is $C_1+C_2+\ol C_1 + \ol C_2$,  after removing it the system is base point free, and that 
it induces a degree 2 morphism $S\to\CP^2$ whose branch divisor is a quartic curve.
Hence again by the diagram \eqref{016} we obtain that
the image under the anticanonical map is a 
scroll of 2-planes over a conic,
and the map is rationally 2 to 1 over the scroll.
Thus $Z$ satisfies all the properties in the proposition.
If $\epsilon_4\circ\epsilon_3\circ\epsilon_2$ blows up one point on each of $C_1, C_2$ 
and  $C_3$, then the line bundle
$K_S^{-1}|_C$ becomes topologically trivial, which again means 
that $Z$ is non-Moishezon.
By symmetries of $S$ we do not need to consider other possibilities, and we have proved the claims 
for the case $k=3$.

If $k=2$, any of the four $\epsilon_i$-s blows up smooth points of the cycle $\epsilon(C)$.
As before write this cycle as $C_1+C_2+\ol C_1+\ol C_2$.
If $\epsilon$ blows up 4 points on $C_1$, then we get $h^0(K_S^{-1})=3$ which is a contradiction.
If $\epsilon$ blows up exactly 3 points on $C_1$, then $S$ is exactly the surface we considered in 
\cite{HonDS4_1} and in particular the system $|2F|$ satisfies the properties in the proposition
by \cite[Proposition 3.2]{HonDS4_1}.
If $\epsilon$ blows up exactly 2 points on $C_1$, then it also blows up exactly 2 points on $C_2$,
which implies that the line bundle $K_S^{-1}|_C$ is topologically trivial.
So again this cannot happen under our Moishezon assumption.
Clearly these cover all possible situations
for the case $k=2$.
\proofend

\begin{remark}
{\em
Propositions \ref{prop:k6}, \ref{prop:k5}, \ref{prop:k4} and \ref{prop:k32} mean that 
 for any Moishezon twistor space
on $4\CP^2$ satisfying $h^0(F)=2$,
the anticanonical image is at least 2-dimensional.
In other words,
the anticanonical system $|2F|$ is always not composed with the pencil $|F|$.
This is the reason why a classification of Moishezon twistor spaces
on $4\CP^2$ can be done by means of the anticanonical system.
}
\end{remark}

\subsection{Summary and conclusions}

In Sections \ref{ss:k6}--\ref{ss:k32}  we have studied Moishezon twistor spaces on $4\CP^2$
which satisfy $h^0(F)=2$.
In order to place those satisfying $h^0(F)>2$ in the same perspective,
we first investigate the anticanonical system of these twistor spaces.

As in Proposition \ref{prop:Kr}, if a twistor space on $4\CP^2$ satisfies 
$h^0(F)>2$, then $Z$ is either a LeBrun twistor space or a Campana-Kreussler twistor space.
For LeBrun twistor spaces,  by putting $n=4$ in \cite[Proposition 2.2 and Corollary 2.3]{HonLB2}, we have
$h^0(2F) =h^0(\CP^1\times\CP^1,\mathscr O(2,2)) = 9$ and the system $|2F|$ is generated by
$|F|$.
Since the system $h^0(F)=4$ and $|F|$ induces a map to a smooth quadric in $\CP^3$,
we obtain the commutative diagram of left
\begin{equation}\label{017}
 \CD
 Z_{\rm{LB}}@>\Phi_{|2F|}>> \CP^8\\
 @V \Phi_{|F|} VV @VV{\iota}V\\
 \CP^3 @>>|\mathscr O(2)|> \CP^9,
 \endCD
\hspace{15mm}
  \CD
 Z_{\rm{CK}}@>\Phi_{|2F|}>> \CP^5\\
 @V \Phi_{|F|} VV @|\\
 \CP^2 @>>|\mathscr O(2)|> \CP^5,
 \endCD
 \end{equation}
where $\iota$ is an embedding as a hyperplane.
In particular, the image $\Phi_{|2F|}(Z)$ is an embedded image  of a smooth quadratic surface in $\CP^3$ under the Veronese image $\CP^3\subset \CP^9$, which is necessarily contained in a  hyperplane in $\CP^9$.

We next see that the situation is also similar  for the Campana-Kreussler twistor spaces.
Let $Z$ be such a twistor space on $4\CP^2$.
Then we have $h^0(F)=3$, and by  \cite[Proposition 1.3 (ii)]{CK98} 
the rational map $Z\to\CP^2$ associated to $|F|$ is surjective.
Also by  \cite[Lemma 1.8]{CK98} we have $h^0(2F)=6$.
These imply that $H^0(2F)$ is generated by $H^0(F)$.
Hence analogously to the left diagram in  \eqref{017} we get the right commutative diagram in \eqref{017}
 Therefore the image by the anticanonical map is the 
Veronese surface $\CP^2$ in $\CP^5$.

Thus in contrast with the case $h^0(F)=2$, for twistor spaces satisfying $h^0(F)>2$,
the anticanonical system $|2F|$ is always generated by $|F|$.

With these remarks, in order to summarize the results in Sections \ref{ss:k6}--\ref{ss:k32}
in a simple form,
we present our results without using the number $k$:
\begin{theorem}\label{thm:smr}
Let $Z$ be a Moishezon twistor spaces on $4\CP^2$.
Then the anticanonical map $\Phi$ of $Z$ satisfies 
exactly one of the following three properties:
\begin{enumerate}
\item
[(I)] $\Phi$ is birational over the image,
\item
[(II)] $\Phi$ is 2 to 1 over the image $\Phi(Z)$,
and the image  is a scroll of 2-planes over a conic,
\item [(III)] the image $\Phi(Z)$ is a rational surface.
\end{enumerate}
\end{theorem}

\proof
The situation (I) happens in Proposition \ref{prop:k6} (i) and (ii) (the case $k=6$),
Proposition \ref{prop:k5} (i) (the case $k=5$) and
Proposition \ref{prop:k4} (i) (the case $k=4$),
and there exists no other Moishezon twistor spaces whose  $\Phi$ is birational,
because the invariant $k$ can catch all Moishezon twistor spaces satisfying
$h^0(F)=2$ by Proposition \ref{prop:cycle}, and for Moishezon twistor spaces with $h^0(F)\neq2$,
$\Phi$ is not birational by the above considerations for the twistor spaces
of LeBrun's and Campana-Kreussler's.
The situation (II) occurs in Proposition \ref{prop:k5} (ii) (the case $k=5$) and
Proposition \ref{prop:k4} (ii) (the case $k=4$) and 
Proposition \ref{prop:k32} (the cases $k=3$ and $k=2$),
and these are all such examples.
The situation (III) happens in Proposition \ref{prop:k4} (iii) (the 
case $k=4$) and also for LeBrun's and Campana-Kreussler's twistor spaces,
and there are no other such examples by our consideration.
\proofend 

\begin{remark}
{\em As above the twistor spaces belonging to (II) 
can be classified into four types according to the value of $k$
(i.e.\,$k=2,3,4,5$).
There is an obvious hierarchy for these twistor spaces; larger $k$ can be obtained from smaller $k$  by taking a limit under smooth deformations.
We remark that, although less trivial, 
 this hierarchy is analogous to the classification in the case of $3\CP^2$ 
obtained by Kreussler-Kurke \cite{KK92}, where 
they classified double solid twistor spaces by configurations of
degree one divisors.
Concerning  the branch divisor
of the anticanonical map $\Phi$ (in the case of $4\CP^2$),
its defining equation is determined  in \cite{HonDS4_1}  for the case $k=2$,
 and \cite{HonDSn1} for the case $k=5$.
Also in these  papers the dimension of the moduli space is computed,
which is $9$-dimensional when $k=2$,
 $4$-dimensional when $k=5$.
 For the cases of $k=3$ and $k=4$, 
 these are yet to be determined.
%
}
\end{remark}
\begin{theorem}\label{thm:birational}
If a Moishezon twistor space $Z$
on $4\CP^2$ is in the case (I) in Theorem \ref{thm:smr}, we have $h^0(F)=9$ or
$h^0(F) = 7$.
If $h^0(F)=7$, the (birational) image $\Phi(Z)\subset\CP^6$
is a complete intersection of three hyperquadrics.
If $h^0(F) = 9$, the (birational) image $\Phi(Z)\subset\CP^8$
is a non-complete intersection of degree $12$.
\end{theorem}

As mentioned after the proof of Proposition \ref{prop:k6},
it is very likely that if $h^0(F)=9$,
$Z$ is a twistor space of Joyce metric.


\begin{theorem}\label{thm:conicbundle}
If a Moishezon twistor space $Z$
on $4\CP^2$ is in the case (III) in Theorem \ref{thm:smr},
$h^0(2F) = 9,6$ or $5$,
and general fibers of $\Phi$ are non-singular rational curves.
Further,
the image rational surface $\Phi(Z)$ can be explicitly described as follows:
(a) If $h^0(F) = 9$, 
it is a smooth hyperplane section of the Veronese embedding of $\CP^3$ 
by $|\mathscr O(2)|$.
(b) If $h^0(F) = 6$, it is a Veronese surface in $\CP^5$.
(c) If $h^0(F) = 5$, it is an intersection of 2 hyperquadrics whose defining equations
can be taken as $x_0x_1=x_2^2$ and $x_3x_4= q(x_0,x_1,x_2)$,
where $q$ is a quadratic polynomial in $x_0,x_1,x_2$.

\end{theorem}
It is likely that the twistor spaces in (c) are exactly the twistor spaces constructed
in \cite{Hon-I}, without assuming the existence of $\mathbb C^*$-action on the twistor spaces:

As a consequence of our classification,
we can determine all Moishezon twistor spaces on $4\CP^2$
which are obtained as a small deformation of LeBrun twistor spaces.
\begin{proposition}\label{prop:LBsmall}
Let $Z$ be a LeBrun twistor space on $4\CP^2$ which is 
generic in the sense that the identity component of
the holomorphic automorphism group is $\mathbb C^*$.
Then if a Moishezon twistor space $Z_t$ is obtained
as a small deformation of $Z$,
then $Z_t$ is a generic LeBrun twistor space,
a Campana-Kreussler twistor space, or one of
the Moishezon twistor spaces studied in \cite{HonDS4_1}.

%
\end{proposition}

The Moishezon twistor spaces studied in \cite{HonDS4_1} can be
characterized by the properties that $h^0(F)=2$ and $k=2$,
where as before $k$ is the half of the number of irreducible components
of the anticanonical cycle $C$.

Note that actually Campana-Kreussler twistor spaces and the ones in \cite{HonDS4_1}
are obtained as a small deformation of generic LeBrun twistor spaces
(\cite[the proof of Theorem 4.2]{Kr97} and \cite[a comment after Proposition 5.4]{HonDS4_1}). 
We just give an outline of a proof of Proposition \ref{prop:LBsmall},
as there is no interesting point in our proof.
By Proposition \ref{prop:cycle} it suffices to show that any Moishezon twistor space 
$Z_t$ on $4\CP^2$
with $h^0(F)=2$ cannot be obtained as 
a small deformation of a LeBrun twistor space $Z$, except the case $k=2$.
For this, let $S_t\in |F|$ be a real irreducible member.
Then $S_t$ contains the cycle $C$ and we know the self-intersection numbers
of each irreducible components of $C$ in $S_t$ as we obtained in Sections 
\ref{ss:k6}--\ref{ss:k32}.
On the other hand for LeBrun twistor spaces
the structure of any real irreducible member $S$ of $|F|$ can be concretely
obtained from LeBrun's construction.
Then we can show that except the case $k=2$, the surface $S_t$ cannot
be obtained as a small deformation of $S$, which shows that 
$Z_t$ cannot be obtained as a small deformation of $Z$.

\vsp
The next result shows that any Moishezon twistor spaces on $4\CP^2$ except
Campana-Kreussler's can be obtained as a small deformation of
the twistor spaces studied in \cite{HonDS4_1},
and
 supports genericity of these twistor spaces. 

\begin{proposition}
If a Moishezon twistor space $Z$ on $4\CP^2$ is not obtained as a 
small deformation of twistor spaces  studied in \cite{HonDS4_1}.
Then $Z$ is a Campana-Kreussler twistor space.
\end{proposition}

We do not give a proof for this, as it can be shown in a similar way to
Proposition \ref{prop:LBsmall}.


%

\end{document}